\documentclass[12pts]{amsart}
\theoremstyle{plain} \newtheorem{Thm}{ }[section]
\title{Syzygies using vector bundles}
\author{Montserrat Teixidor i Bigas}
\address{Mathematics Department, Tufts University, Medford
MA02155, USA}

\begin{document}

\begin{abstract}
This paper studies syzygies of curves that have been embedded in
projective space by line bundles of large degree. The proofs take
advantage of the relationship between syzygies and spaces of
section of vector bundles associated to the given line bundles.
\end{abstract}

\maketitle

\begin{section}{Introduction}
Let $C$ be a projective irreducible non-singular curve of genus
$g$ defined over a field of characteristic zero. Let $L$ be a line
bundle of degree at least $2g+1$ on $C$ and consider the immersion
of $C$ in projective space associated to the complete linear
series $|L|$. Define the rings
$$R=\oplus _{n=0}^{\infty}H^0({\bf P}, {\mathcal O}(n)),
S=\oplus _{n=0}^{\infty}H^0(C, L^n)$$
Then, $R$ is a graded $S$
module and admits a minimal graded free resolution.
$$0\rightarrow F_t\rightarrow ...\rightarrow F_1\rightarrow F_0\rightarrow 0$$
The vector space of  syzygies $K_{p,q}$ is defined as the piece of
degree $q$ of $F_p$.

In the case of a curve, one easily sees (cf.for example \cite{Ei}
or \ref{p,0,p,3}) that $K_{p,q}=0$ if $q\ge 3$ or $q=0, p\not= 0$.
Hence, one is left with only two strands corresponding to $q=1,2$.
Moreover, these two strands are related (cf. for example
\ref{p+1,1-p,2}), so in fact it suffices to compute one of them.

For curves of high degree, one expects that $K_{p,1}=0$ if $p$ is
large. On the other hand, if the curve has special linear series
$L_1,L_2$ of dimensions $r_1,r_2$ and $L=L_1\otimes L_2$, then
from the converse of Green's conjecture (see \cite{GL1}) ,
$K_{r_1+r_2-1,1}\not= 0$.

The question then arises as to whether the non-vanishing of the
$K_{p,1}$ can be explained by the existence of certain linear
series on the curve. A conjecture of Green and Lazarsfeld (see
\cite{GL2} , conjecture 3.7) states that this is the case : Let
$L$ be a line bundle of sufficiently high degree giving rise to a
linear series of dimension $r$, then $K_{r-d,1}(L)=0$ unless $C$
is $d$-gonal.

This conjecture has been proved recently for the generic curve
(\cite{AV}, \cite{A}) or equivalently, for the curves of largest
possible gonality. In this paper we deal with the other end of the
spectrum, namely the curves with small gonality. In \ref{r-2,1},
we provide a new proof of Green's $K_{p,1}$ Theorem (see \cite{G}
Th 3.c.1, p.151) in the case of curves. In \ref{r-3,1}, we reprove
Ehbauer's Theorem (\cite{Eh} Th 1.4, p.146 ) . While the results
contained here are not new, we do provide new proofs and a
simpler, unified presentation.
\end{section}

\begin{section}{Preliminaries}

It is standard, that the $K_{p,q}$ can be computed by taking
homology in the middle term of the Koszul complex
$$\wedge ^{p+1} H^0(L)\otimes H^0(C,L^{q-1}) \rightarrow
\wedge ^{p} H^0(L)\otimes H^0(C,L^{q} )\rightarrow \wedge ^{p-1}
H^0(L)\otimes H^0(C,L^{q+1} )\rightarrow...$$ Define now a  vector
bundle  $E_L$ by using the exact sequence
$$0\rightarrow E_L^*\rightarrow H^0(L)\otimes {\mathcal O}_C
\rightarrow L \rightarrow 0.$$ In order to compare this with the
Koszul complex, we consider the wedge powers of the sequence above

$$(*)0\rightarrow \wedge ^{p}E_L^*\otimes L^{q}\rightarrow \wedge ^{p}H^0(L)\otimes
L^{q} \rightarrow \wedge^{p-1} E_L^*\otimes L ^{q+1}\rightarrow
0$$ and the analogous sequences obtained by replacing $p,q$ by
$p+1,q-1$ and $p-1,q+1$. The map in the Koszul sequence
$$...  \wedge^{p+1}H^0(L)\otimes H^0(L^{q-1})\rightarrow \wedge^{p}H^0(L)\otimes H^0(L^{q})
...,$$ can be factored through $H^0(\wedge^p E^*\otimes L ^q)$.
The two factoring  maps
$$\wedge^{p+1}H^0(L)\otimes H^0(L^{q-1})\rightarrow H^0(\wedge^p E_L^*\otimes L
^q)$$ and
$$H^0(\wedge^p E_L^*\otimes L ^q)\rightarrow \wedge^{p}H^0(L)\otimes
H^0(L^{q})$$ are obtained by taking homology on the sequences
above. One obtains then the following presentations for the
syzygies

\begin{Thm} \label{p,q} {\bf Lemma.}
$$K_{p,q}={H^0(\wedge ^p E_L^*\otimes L^q)\over Im [\wedge
^{p+1}H^0(L)\otimes H^0(L^{q-1})]}=Im[H^0(\wedge ^p E_L^*\otimes
L^q)\rightarrow H^1(\wedge ^{p+1} E_L^*\otimes L^{q-1})]$$

\end{Thm}

{\bf Notation}

In all that follows, we will assume that $d\ge 2g+1$ and that  $L$
is a line bundle of degree $d$ on $C$. We define $r$ by
$h^0(L)=d+1-g=r+1$. We have a natural immersion $C\rightarrow {\bf
P}^r$. We will use freely the identifications given in \ref{p,q}
for various values of $p,q$.

We shall write $E$ instead of $E_L$ when $L$ is clear from the
context.

The first statement in the Lemma below was proved by David Butler
in \cite{B} Theorem 1.2. The second follows from the first because
in characteristic zero the wedge powers of a stable bundle are
semistable

\begin{Thm}\label{Butler}{\bf Lemma.} Under the above conditions
on degree, the vector bundle $E_L$ is stable. Hence, in
characteristic zero, its wedge powers are semistable.
\end{Thm}

\end{section}

\begin{section}{Syzygies of curves of large degree}

\begin{Thm}\label{p,0,p,3}{\bf Proposition} $K_{p,0}=0, p\not= 0,\
dim K_{0,0}=1$, $K_{p.q}=0$ if  $q\ge 3$.
\end{Thm}
\begin{proof}
The exact sequence (*) corresponding to  $p+1,0$ becomes
$$0\rightarrow \wedge ^{p+1}E^*\otimes L^{-1}\rightarrow \wedge ^{p+1}H^0(L)\otimes
L^{-1} \rightarrow \wedge^p E^*\rightarrow 0.$$ If $p> 0$,
$h^0(\wedge ^pE^*)=0$ as this is a semistable sheaf of negative
degree. Hence, from the first description of $K_{p,0}$ in
\ref{p,q}, $K_{p,0}=0$.

 If $p=0$,  $H^0(\wedge^pE^*)=H^0(\mathcal O)$ has dimension one while $H^0(L)\otimes
 H^0(L^{-1})=0$. Hence, $dim K_{0,0}=1$.

Assume now $p\not= 0,\ q\ge 3$. As $rank E=d-g=r$, if $p+1 >r$,
$\wedge^{p+1}E=0$. If $p+1\le r$,  the slope
$$\mu (\wedge ^{p+1}E^*\otimes L^{q-1})=(q-1)d-\frac {(p+1)d}
{r}=d(q-1-\frac {p+1}{ r})\ge d>2g-2.$$ Hence, $h^1(\wedge
^{p+1}E^*\otimes L^{q-1})=0$. Then, from the second description of
$K_{p,q}$ in \ref{p,q}, the result follows.
\end{proof}

\begin{Thm}\label{p+1,1-p,2}{\bf Proposition}
$dimK_{p.2}-dimK_{p+1,1}$ depends only on $d$ and not on $L$ (and
can be explicitly calculated).
\end{Thm}
\begin{proof}
Taking $q=1$ and $p+1$ instead of $p$ in the exact sequence (*),
we get
$$0\rightarrow \wedge ^{p+2}E^*\rightarrow \wedge ^{p+2}H^0(L)\otimes
{\mathcal O}\rightarrow \wedge^{p+1} E^*\otimes L \rightarrow 0$$
 As $\wedge ^{p+2}E^*$ is semistable of negative degree, its space of sections is zero.
Hence, from the first description of $K_{p+1,1}$ in \ref{p,q}

 $$dim K_{p+1,1}=h^0(\wedge ^{p+1}E^*\otimes L)-{r+1\choose p+2}.$$
On the other hand, from
$$0\rightarrow \wedge ^{p+1}E_L^*\otimes L\rightarrow \wedge ^{p+1}H^0(L)\otimes
L\rightarrow \wedge^{p} E^*\otimes L^2 \rightarrow 0,$$

 $$dim K_{p,2}=h^0(\wedge ^pE^*\otimes L^2)-(r+1){r+1\choose
 p+1}+h^0(\wedge ^{p+1}E^*\otimes L).$$
 Therefore, $$dim K_{p,2}-dim K_{p+1,1}=h^0(\wedge ^pE^*\otimes L^2)+{r+1\choose p+2}-(r+1){r+1\choose
 p+1}.$$
 As $\wedge ^pE^*\otimes L^2$ has slope greater than $2g-2$, the dimension of its space of sections is
 independent of $L$ for a fixed $d$,
 namely $h^0(\wedge ^pE^*\otimes L^2)={r\choose p}(2d-{pd\over
 r}+1-g)$. Then $$K_{p,2}-K_{p+1,1}={r\choose p}(2d-{pd\over
 r}+1-g)+{r+1\choose p+2}-(r+1){r+1\choose
 p+1}.$$

 \end{proof}

\begin{Thm}\label{r-1,2}{\bf Proposition.}
$dimK_{r-1,2}=g$
\end{Thm}

\begin{proof} (see also \cite{Ei} Prop 8.6)
As $h^0(L)=r+1$, $\wedge ^{r+1}H^0(L)$ is isomorphic to the base
field and $\wedge ^rH^0(L)$ is naturally isomorphic to
$(H^0(L))^*$. As $rk E=r$ and $\wedge ^{r}E=L$, then
$\wedge^{r-1}E^*=E\otimes L^{-1}$. Using these isomorphisms, the
exact sequence
$$0\rightarrow \wedge^{r}E^*\otimes L\rightarrow \wedge ^{r}
H^0(L)\otimes L \rightarrow \wedge ^{r-1}E^*\otimes L^2\rightarrow
0$$ can be written as
$$0\rightarrow {\mathcal O}\rightarrow H^0(L)^*\otimes
L\rightarrow E\otimes L\rightarrow 0.$$ Taking homology, one
obtains
$$H^0(L)^*\otimes H^0(L)\rightarrow H^0(E\otimes L)\rightarrow
H^1({\mathcal O})\rightarrow 0$$ Hence, $dim
K_{r-1,2}=h^1({\mathcal O})=g$
\end{proof}

{\bf Remark.} From the second interpretation in \ref{p,q},
$K_{p,q}=0$ if $p\ge r$. Then,  the result also follows from
\ref{p+1,1-p,2}.

\begin{Thm} \label{r-1,1} {\bf Proposition} For $C$ a non-rational
curve,  $K_{r-1,1}=0$.
\end{Thm}

\begin{proof}
Again, we use the isomorphisms $\wedge^kE^*=\wedge^{r-k}E\otimes
L^{-1}, \wedge^kH^0(L)=\wedge^{r+1-k}H^0(L)^*$. The exact sequence
$$0\rightarrow \wedge ^{r}E^*\rightarrow \wedge
^{r}H^0(L)\otimes {\mathcal O}\rightarrow \wedge ^{r-1}E^*\otimes
L\rightarrow 0$$ becomes
$$0\rightarrow L^*\rightarrow H^0(L)^*\otimes {\mathcal
O}\rightarrow E\rightarrow 0.$$ Taking homology, this gives
$$0\rightarrow H^0(L)^*\rightarrow H^0(E)\rightarrow H^1(L^*)$$
Hence, the statement is equivalent to the surjectivity of the map
$H^0(L)^*\rightarrow H^0(E)$. By the injectivity of this map, it
suffices to prove that $h^0(E)=r+1$. This fact is proved as
follows (see \cite{GL3}, diagram 2.1): Let $P_1,...P_{r-1}$ be
generic points, $D=P_1+...+P_{r-1}$. Then, $h^0(L(-D))=2$.
Consider the diagram
$$\begin{matrix} 0&\rightarrow &L^{-1}(D)&\rightarrow &
H^0(L(-D)) \otimes  {\mathcal O}&\rightarrow & L(-D)&
\rightarrow& 0\\
 & &\downarrow &  &\downarrow&  & \downarrow& &\  \\
 0&\rightarrow &E^*&\rightarrow &H^0(L)\otimes {\mathcal
O}&\rightarrow &L&\rightarrow&0\\
 & &\downarrow &  &\downarrow& &  \downarrow& &\  \\
0&\rightarrow &\bar E^*&\rightarrow &{H^0(L)\over
H^0(L(-D))}\otimes {\mathcal O}&\rightarrow &L_D&\rightarrow&0.\\
\end{matrix}$$
Here $\bar E^*$ is defined as the cokernel of the left vertical
sequence. It can be shown then that it is also the kernel of the
bottom row. From this description, it follows that $\bar
E^*=\oplus_{i=1}^{r-1}{\mathcal O}(-P_i)$.  We then get an exact
sequence
$$ 0\rightarrow \oplus_{i=1}^{r-1} {\mathcal O}(P_i)\rightarrow E\rightarrow
L(-D)\rightarrow 0.$$ Therefore,
$$h^0(E)\le \sum h^0({\mathcal
O}(P_i))+h^0(L(-D))=r-1+2=r+1$$ as required.
\end{proof}

\begin{Thm}\label{p,2} {\bf Proposition.}
 Let $C$ be a curve immersed by a line bundle $L$ of degree
$2g+1+k$. Then, if $p\le k, K_{p,2}=0$.
\end{Thm}
\begin{proof} (See \cite{L}, section 1)
From the second interpretation of $K_{p,2}$ in \ref{p,q},
 it suffices to show that
$h^1(\wedge^{p+1}E^*\otimes L)=0$.

From the proof of \ref{r-1,1}, we have an exact sequence
$$0\rightarrow L^{-1}(D)\rightarrow E^*\rightarrow \oplus
_{i=1}^{r-1}{\mathcal O}(-P_i)\rightarrow 0$$ where $D= \oplus
_{i=1}^{r-1}P_i$ is a divisor made of generic points. Taking wedge
powers, we obtain exact sequences
$$0\rightarrow \oplus L^{-1}(D-P_{i_1}-...-P_{i_p})\rightarrow
\wedge ^{p+1}E^*\rightarrow \oplus {\mathcal
O}(-P_{i_1}-...-P_{i_{p+1}})\rightarrow 0.$$ Tensoring with $L$,
this gives rise to
$$0\rightarrow \oplus {\mathcal O}(D-P_{i_1}-...-P_{i_p})\rightarrow
\wedge ^{p+1}E^*\otimes L\rightarrow \oplus
L(-P_{i_1}-...-P_{i_{p+1}})\rightarrow 0.$$

Given a choice of points $P_{i_1}...P_{i_{p+1}}$ denote by
$P_{j_1}...P_{j_{r-1-p}}$ the complementary set of points in $D$.
Then,
$h^1(D(-P_{i_1}-...-P_{i_{p+1}}))=h^1(P_{j_1}+...+P_{j_{r-1-p}})=0$
if and only if (by the genericity of the points) $r-1-p\ge g$.
This condition can be written as $p\le r-1-g=d-2g-1=k$ and this is
satisfied by assumption. Moreover, $deg K\otimes
L^{-1}(P_{i_1}+...+P_{i_{p+1}})=p-k-2<0$. Then, $0=h^0(K\otimes
L^{-1}(P_{i_1}+...+P_{i_{p+1}}))=h^1(L(-P_{i_1}-...-P_{i_{p+1}}))$.
Hence, from the exact sequence above, $h^1(\wedge ^{p+1}E^*\otimes
L)=0$ if $p\le k$.
\end{proof}

We have seen in \ref{r-1,1} that $K_{r-1,1}=0$. We want to see
that $K_{p,1}=0$ if $p$ is large enough and the curve is
sufficiently general. The following Lemma will be useful:

\begin{Thm}\label{lemap,1}
{\bf Lemma} $K_{p,1}=0$ if and only if
$h^0(\wedge^{r-p}E)={r+1\choose r-p}$.
\end{Thm}
\begin{proof}
From \ref{p,q} , $K_{p,1}=0$ if and only if the natural map
$\wedge^{p+1}H^0(L)\rightarrow H^0(\wedge ^pE^*\otimes L) $ is
surjective.
 The kernel of this map is $H^0(\wedge^{p+1}E^*)$ and this space
 of sections is zero because $\wedge^{p+1}E^*$ is a semistable
 vector bundle of negative degree. Therefore, $K_{p,1}=0$ if and only if
 $h^0(\wedge ^pE^*\otimes L)={r+1\choose p+1}$.

 As $rank E=r$ and $\wedge ^rE=L$, one has an isomorphism of $\wedge ^pE^*\otimes
 L$ with $\wedge ^{r-p}E$. Moreover, ${r+1\choose p+1}={r+1\choose
 r-p}$ and the result follows.

\end{proof}

 The proof of the next two Propositions is inspired in the paper
 of Claire Voisin \cite{V}

\begin{Thm} \label{r-2,1} {\bf Proposition.} Let $L$ be a line bundle of degree at
least $2g+1$ on a curve of genus at least four. Then $K_{r-2,1}=0$
unless $C$ is trigonal and $L=K(g^1_3)$ or $C$ is hyperelliptic.
\end{Thm}
\begin{proof}
We shall see at the end that the assumptions on $L$ imply that
there exists a divisor $D=P_1+...+P_r$ satisfying the following
conditions

a) $h^0(L(-D))=h^0(L(-(D-P_i))=2$ for every $i=1...r$

b) $h^0(P_{i_1}+P_{i_2}+P_{i_3})=1$ for all triples of points in
the support of $D$.

 We now show that these conditions imply that $h^0(\wedge^2(E))={r+1\choose 2}$
 and therefore $K_{r-2,1}=0$ by \ref{lemap,1}.

Consider the following exact diagram
$$\begin{matrix} 0&\rightarrow &(E_{L(-D)})^*&\rightarrow &
H^0(L(-D))\otimes {\mathcal O}&\rightarrow &L(-D)&\rightarrow&0\\
 & &\downarrow & &\downarrow& & \downarrow& & \\
 0&\rightarrow &E^*&\rightarrow &H^0(L)\otimes {\mathcal
O}&\rightarrow &L&\rightarrow&0\\
 & &\downarrow & &\downarrow& & \downarrow& & \\
0&\rightarrow & \bar E^*&\rightarrow&V\otimes {\mathcal
O}&\rightarrow &L_D&\rightarrow&0\\
\end{matrix}$$
where $V={H^0(L)\over H^0(L(-D))}$.

 As
$h^0(L(-D))=2$,  $(E_{L(-D)})^*=L^{-1}(D)$.  The condition
 $h^0(L(-(D-P_i)))=h^0(L)-deg(D-P_i)$
implies that $h^0(L(-(D-P_i-P_j)))=h^0(L(-(D-P_i)))+1=3$. Take a
section $s$ of $H^0(L(-(D-P_1-P_2)))-H^0(L(-(D-P_i))$. The image
of $s$ in the quotient $V$ is non-zero and $dim{V\over<s>}=r-2$.
We then have
$$\begin{matrix} 0&\rightarrow &{\mathcal O}(-P_1-P_2)&\rightarrow &
<s>\otimes {\mathcal O}&\rightarrow &L_{P_1+P_2}&\rightarrow&0\\
 & &\downarrow & &\downarrow& & \downarrow& & \\
 0&\rightarrow &\bar E^*&\rightarrow &V\otimes {\mathcal
O}&\rightarrow &L_D&\rightarrow&0\\
 & &\downarrow & &\downarrow& & \downarrow& & \\
0&\rightarrow & A&\rightarrow&{V\over <s>}\otimes {\mathcal
O}&\rightarrow &L_{P_3+...+P_r}&\rightarrow&0.\\
\end{matrix}$$
From the diagram, $A={\mathcal O}(-P_3)\oplus...\oplus {\mathcal
O}(-P_r)$. We then get exact sequences
$$0\rightarrow \oplus _{i=3}^r {\mathcal O}(P_i)\rightarrow \bar
E\rightarrow {\mathcal O}(P_1+P_2)\rightarrow 0$$
$$0\rightarrow  \bar E\rightarrow  E\rightarrow
L(-D)\rightarrow 0.$$ Taking the second wedge power of the latter,
we obtain
$$0\rightarrow  \wedge ^2\bar E\rightarrow  \wedge ^2E\rightarrow
\bar E\otimes L(-D)\rightarrow 0.$$ From the first exact sequence,
we get
$$0\rightarrow \oplus _{i=3}^rL(-(D-P_i))\rightarrow \bar
E\otimes L(-D)\rightarrow L(-D+P_1+P_2)\rightarrow 0$$
$$0\rightarrow \oplus _{3\le i_1<i_2\le r} {\mathcal O}(P_{i_1}+P_{i_2})\rightarrow \wedge ^2\bar
E\rightarrow \oplus_{i=3}^r {\mathcal O}(P_1+P_2+P_i)\rightarrow
0.$$

Hence,
$$h^0(\wedge ^2E)\le h^0(\wedge ^2\bar E)+h^0(\bar E\otimes
L(-D))\le \sum  _{3\le i_1<i_2\le r} h^0({\mathcal
O}(P_{i_1}+P_{i_2}))+$$ $$\sum _{i=3}^r h^0({\mathcal
O}(P_1+P_2+P_i))+\sum
_{i=3}^rh^0(L(-(D-P_i)))+h^0(L(-D+P_1+P_2))=$$
$${r-2\choose 2}+r-2+2(r-2)+3={r-2\choose 2}{3\choose 0}+{r-2\choose 1}{3\choose 1}+{r-2\choose 0}{3\choose
2}= {r+1\choose 2}.$$
\bigskip

It remains to show that a divisor $D$ exists satisfying conditions
a),b) above. Let $D_{g-2}$ be a generic effective divisor of
degree $g-2$. Let $D$ be an effective divisor of the complete
linear series $|L\otimes K^{-1}(D_{g-2})|$. This exists as
$deg(L\otimes K^{-1}(D_{g-2}))=d-g\ge g+1$ and therefore $L\otimes
K^{-1}(D_{g-2})$ has sections. Moreover,
$h^0(L(-D))=h^0(K(-D_{g-2}))=2$. We now show that for every point
$P_i$ in the support of $D$, $h^0(L(-(D-P_i)))=2$. As
$L(-(D-P_i))=K(-D_{g-2}+P_i)$, this is equivalent to
$h^0(D_{g-2}-P_i)=0$. From our choice, $h^0(D_{g-2})=1$. Hence we
are asking that $D, \ D_{g-2}$ have disjoint supports.

  Such a $D$ will exist if $h^0(L\otimes
K^{-1}(D_{g-2}))>h^0(L\otimes K^{-1}(D_{g-2}-Q))$ for each $Q$ in
the support of $D_{g-2}$. This condition is equivalent to
$h^0(K^2\otimes L^{-1}(-D_{g-2}))=h^0(K^2\otimes
L^{-1}(-(D_{g-2}-Q))$. By the genericity of $D_{g-2}$ this happens
precisely when
 $h^0(K^2\otimes L^{-1})\le g-3$.  Assume on the contrary that
 $h^0(K^2\otimes L^{-1})\ge g-2$. We shall prove that this leads
 to one of the situations we excluded for $L$. We obtain
 $$g-2\le h^0(K^2\otimes L^{-1})=h^0(L\otimes K^{-1})+3g-3-d\le h^0(L\otimes K^{-1})+g-4 .$$
Hence, $h^0(L\otimes K^{-1})\ge 2$ . It follows that $K^2\otimes
L^{-1}$ contributes to the Clifford index of $C$ as $g\ge 4$ and
$$0\le Cliff(K^2\otimes L^{-1})=4g-4-d-2(h^0(K^2\otimes
L^{-1})-1)\le 4g-4-2g-1-2(g-3)$$ Therefore, the only two
possibilities are $degL=2g+1,\ h^0(K^2\otimes L^{-1})=g-2$ and
$degL=2g+2,\ h^0(K^2\otimes L^{-1})=g-2$. In the first case,
$L\otimes K^{-1}$ gives a $g^1_3$, in the second it gives a
$g^2_4$. As $g^2_4$ appear only on hyperelliptic curves and
$g^2_4=2g^1_2$, both cases have been excluded. This takes care of
condition a).

Condition b) is obviously satisfied if $C$ is not trigonal. If $C$
is trigonal, we can choose a $D$ satisfying b) so long as the
linear system $|D|$ is not composed with  a $g^1_3$. If we assume
the curve non-hyperelliptic, the $g^1_3$ is unique for $g\ge 5$
(there may be two of them for genus $g=4$). If $|D|$ is composed
with the $g^1_3$,  then $L\otimes K^{-1}(D_{g-2})=ag^1_3, a\ge
{g+2\over 3}$. As the divisor $D_{g-2}$ depends on $g-2$
parameters, this is impossible.

\end{proof}

\bigskip

\begin{Thm}\label{r-3,1}  {\bf Proposition}
Let $L$ be a line bundle of degree $d\ge 2g+3$ on a curve $C$ that
is not trigonal of genus at least seven. Then $K_{r-3,1}=0$
\end{Thm}
\begin{proof}
We shall check in a moment that the given conditions imply that
there is an effective divisor of degree $r-1$  $D=P_1+...+P_{r-1}$
such that

 a) $h^0(L(-D))=3= h^0(L(-D+P_i)),\ i=1...r-1$.

 b) $h^0(P_{i_1}+P_{i_2}+P_{i_3}+P_{i_4})=1$ for all quadruples of
 points on the support of $D$.

 c) The maps $H^0(L(-D))\otimes H^0(K(-D'))\rightarrow H^0(K\otimes
 L(-D-D'))$ are onto for every divisor $D'$ of degree at most
 three contained in the support of $D$.

Using conditions a)-c), we now show that
$h^0(\wedge^3E)={r+1\choose 3}$. From \ref{lemap,1}, this suffices
to prove the result.

 Consider the diagram
$$\begin{matrix} 0&\rightarrow &(E_{L(-D)})^*&\rightarrow &
H^0(L(-D))\otimes {\mathcal O}&\rightarrow &L(-D)&\rightarrow&0\\
 & &\downarrow & &\downarrow& & \downarrow& & \\
 0&\rightarrow &E^*&\rightarrow &H^0(L)\otimes {\mathcal
O}&\rightarrow &L&\rightarrow&0\\
 & &\downarrow & &\downarrow& & \downarrow& & \\
0&\rightarrow & \bar E^*&\rightarrow&V\otimes {\mathcal
O}&\rightarrow &L_D&\rightarrow&0\\
\end{matrix}$$
where $V=\frac {H^0(L)}{H^0(L(-D))}$ is a vector space of
dimension $r-2$. We then have the exact sequence
$$0\rightarrow \bar E \rightarrow E\rightarrow E_{L(-D)}
\rightarrow 0$$ and $E_{L(-D)}$ is a vector bundle of rank two.

 We
obtain the bound
$$h^0(\wedge^3E)\le h^0(\wedge ^3\bar E)+h^0(\wedge ^2\bar E
\otimes E_{L(-D)})+h^0(\bar E \otimes \wedge ^2E_{L(-D)}).$$

 As
$h^0(L(-D))=h^0(L)-deg D+1$ and $h^0(L(-(D-P_i)))=h^0(L(-D))$, it
follows that $h^0(L(-(D-P_i-P_j)))=h^0(L(-(D-P_i)))+1=4$. Take a
section $s$ of $H^0(L(-(D-P_1-P_2)))-H^0(L(-(D-P_i))$. The image
of $s$ in the quotient $V$ is non-zero and $dim{V\over<s>}=r-3$ We
 have
$$\begin{matrix} 0&\rightarrow &{\mathcal O}(-P_1-P_2)&\rightarrow &
<s>\otimes {\mathcal O}&\rightarrow &L_{P_1+P_2}&\rightarrow&0\\
 & &\downarrow & &\downarrow& & \downarrow& & \\
 0&\rightarrow &\bar E^*&\rightarrow &V\otimes {\mathcal
O}&\rightarrow &L_D&\rightarrow&0\\
 & &\downarrow & &\downarrow& & \downarrow& & \\
0&\rightarrow & \oplus _{i=3}^{r-1}{\mathcal
O}(-P_i)&\rightarrow&{V\over <s>}\otimes {\mathcal O}&\rightarrow
&L_{P_3+...+P_{r-1}}&\rightarrow&0.\\
\end{matrix}.$$

Therefore, the following sequences are exact
$$0\rightarrow \oplus _{i=3}^{r-1} {\mathcal O}(P_i)\rightarrow \bar
E\rightarrow {\mathcal O}(P_1+P_2)\rightarrow 0$$
$$0\rightarrow L^{-1}(D)\rightarrow H^0(L(-D))^*\otimes {\mathcal
O}\rightarrow E_{L(-D)}\rightarrow 0.$$

Taking wedge powers, we get
$$0\rightarrow \oplus {\mathcal O}(P_{i_1}+P_{i_2}+P_{i_3})\rightarrow
\wedge ^3\bar E\rightarrow \oplus {\mathcal
O}(P_{i_1}+P_{i_2}+P_1+P_2)\rightarrow 0.$$ Using condition b)
$$h^0(\wedge^3\bar E)\le {r-3\choose 3}+{r-3\choose 2}.$$

Also,  as $\wedge ^2E_{L(-D)}=L(-D)$,
$$h^0(\bar E \otimes \wedge^2E_{L(-D)})\le
\sum_{i=3}^{r-1}h^0(L(-D+P_i))+h^0(L(-D+P_1+P_2))=3(r-3)+4.$$

From
$$0\rightarrow \oplus {\mathcal O}(P_{i_1}+P_{i_2})\rightarrow
\wedge ^2\bar E\rightarrow \oplus {\mathcal
O}(P_{i}+P_1+P_2)\rightarrow 0$$, one has $$h^0(\wedge ^2\bar
E\otimes E_{L(-D)})\le
\sum_{i_1,i_2}h^0(E_{L(-D)}(P_{i_1}+P_{i_2}))+\sum_{i}h^0(E_{L(-D)}(P_1+P_2+P_i)).$$
In order to compute the second of these numbers, we use the exact
sequence

$$0\rightarrow L^{-1}(D+P_1+P_2+P_i)\rightarrow H^0(L(-D))^*\otimes {\mathcal
O}(P_1+P_2+P_i)\rightarrow E_{L(-D)}(P_1+P_2+P_i)\rightarrow 0$$

 The map $H^1(L^{-1}(D+P_1+P_2+P_i))\rightarrow H^0(L(-D))^*\otimes
H^1({\mathcal O}(P_1+P_2+P_i))$ is dual of the  map
$$H^0(K\otimes L(-D-P_1-P_2-P_i))\leftarrow H^0(L(-D))\otimes
H^0(K(-P_1-P_2-P_i))$$ which is surjective by assumption c). As
$h^0(L^{-1}(D+P_1+P_2+P_i))=0$, we obtain
$$h^0(E_{L(-D)}(P_1+P_2+P_i))\le h^0(L(-D))h^0({\mathcal
O}(P_1+P_2+P_i))=3.$$ Similarly

$$h^0(E_{L(-D)}\otimes {\mathcal O}(P_{i_1}+P_{i_2}))\le h^0(L(-D))h^0({\mathcal
O}(P_{i_1}+P_{i_2}))=3.$$ Hence,

$$h^0(\wedge ^2\bar E\otimes E_{L(-D)})\le 3{r-3\choose
2}+3(r-3).$$ We then deduce that $$h^0(\wedge ^3E_L)\le
{r-3\choose 3}+{r-3\choose 2}+4+3(r-3)+3{r-3\choose 2}+3(r-3)=$$
$${4\choose 0}{r-3\choose 3}+{4\choose 1}{r-3\choose 2}+ {4\choose
2}{r-3\choose 1}+{4\choose 3}{r-3\choose 0}={r+1\choose 3}.$$
\bigskip

It remains to show that a divisor $D$ exists satisfying conditions
a),b),c) above. Let $D_{g-3}$ be a generic effective divisor of
degree $g-3$. Choose as  $D$ a generic effective divisor of the
complete linear series $|L\otimes K^{-1}(D_{g-3})|$. This is
possible as $deg(L\otimes K^{-1}(D_{g-3}))=d-g-1\ge g+2$ .
Moreover, $h^0(L(-D))=h^0(K(-D_{g-3}))=3$. We now show that we can
choose $D$ so that for every point $P_i$ in the support of $D$,
$h^0(L(-(D-P_i)))=3$. As $L(-(D-P_i))=K(-D_{g-3}+P_i)$, this is
equivalent to $h^0(D_{g-3}-P_i)=0$. From the genericity of
$D_{g-3}$, $h^0(D_{g-3})=1$. Hence we are asking that $D, \
D_{g-3}$ have disjoint supports.

  Such a $D$ will exist if $h^0(L\otimes
K^{-1}(D_{g-3}))>h^0(L\otimes K^{-1}(D_{g-3}-Q))$ for each $Q$ in
the support of $D_{g-3}$. By Serre duality, this translates into
$h^0(K^2\otimes L^{-1}(-D_{g-3}))=h^0(K^2\otimes
L^{-1}(-(D_{g-3}-Q))$ . By the genericity of $D_{g-3}$, this is
equivalent to
 $h^0(K^2\otimes L^{-1})\le g-4$. Assume the opposite, namely
 $h^0(K^2\otimes L^{-1})\ge g-3$. Then,
$$h^0(L\otimes K^{-1})\ge d-2g+2+1-g+g-3=d-2g\ge 2.$$
 It follows that
$K^2\otimes L^{-1}$ contributes to the Clifford index of $C$ and
$$0\le Cliff(K^2\otimes L^{-1})=4g-4-d-2(h^0(K^2\otimes
L^{-1})-1)\le 2.$$ Hence,
$$d\le 4g-4-2(h^0(K^2\otimes L^{-1})-1)\le 2g+4.$$
Then,
$$h^0(L\otimes K^{-1})=d+3-3g+h^0(K^2\otimes L^{-1})\ge d-2g.$$
As we are assuming $degL\ge 2g+3$, the only possibilities are then
 $degL= 2g+3 ,\ h^0(L\otimes K^{-1})\ge 3$ and $degL= 2g+4 ,\
h^0(L\otimes K^{-1})\ge 4$.

Therefore, $L=K(g^2_5)$ which contradicts the bound on the genus
of $C$ or $L=K(g^3_6)$ and then, the curve is hyperelliptic which
contradicts the assumption. This takes care of condition a).

Condition b) is obviously satisfied if $C$ is not $4$-gonal. Note
that $dim W^1_4\le 1$ for $C$ non-hyperelliptic and equality holds
only for $C$ trigonal or bielliptic. If $C$ is $4$-gonal, we can
choose a $D$ satisfying b) so long as $L\otimes K^{-1}(D_{g-3})$
is not composed with a $g^1_4$.  Otherwise, for $C$ non-trigonal
or bielliptic, this would imply that $L\otimes
K^{-1}(D_{g-3})=kg^1_4$ . This is incompatible with the fact that
$L\otimes K^{-1}(D_{g-3})$ moves in a $(g-3)$-dimensional family
of line bundles. For $C$ bielliptic, $W^1_4$ is the pull-back of
the linear series of degree two. Hence, if the $D_i$ are divisors
of various $g^1_4$, $h^0(D_1+...+D_k)\ge 2k>k+1$ if $k>1$. This
contradicts the statement that $L\otimes K^{-1}(D_{g-3})$ is
composed with an involution.

\bigskip
Before proving c), we show that $|D|$ has no fixed points and is
not composed with an involution. Assume that $|D|$ had fixed
points for generic $D_{g-3}$. Then for every divisor $D_{g-3}$
there exists a $P_{D_{g-3}}$ such that
$$h^0(L\otimes
K^{-1}(D_{g-3}))=h^0(L\otimes K^{-1}(D_{g-3}-P_{D_{g-3}})).$$
Equivalently

$$h^0(K^2\otimes L^{-1}(-D_{g-3}))=h^0(K^2\otimes L^{-1}(P_{D_{g-3}}-D_{g-3}))-1.$$
This implies in particular that $K^2\otimes
L^{-1}(P_{D_{g-3}}-D_{g-3})$ is effective. As this is a family of
dimension at least $g-4$ of line bundles, if they are effective
they must have degree at least $g-4$. Hence,
$$4g-4-d+1-g+3\ge g-4$$
Therefore $d\le 2g+4$.

 Let us show that the condition $h^0(K^2\otimes
L^{-1}(P_{D_{g-3}}-D_{g-3}))\ge 1$ implies $h^0(K^2\otimes
L^{-1})\ge g-3$: If $P_{D_{g-3}}=P$ is fixed independent of
$D_{g-3}$, then $h^0(K^2\otimes L^{-1}(P))\ge g-2$. Hence,
$h^0(K^2\otimes L^{-1})\ge g-3$. If $P_{D_{g-3}}$ is generic, then
$h^0(K^2\otimes L^{-1}(P_{D_{g-3}}-D_{g-3}))=h^0(K^2\otimes
L^{-1}(-D_{g-3}))\ge 1$. Moreover, for a fixed $P$, $D_{g-3}$
moves in a variety of dimension at least $g-4$. Hence,
$h^0(K^2\otimes L^{-1})\ge g-3$.

If $d=2g+4$, from $h^0(K^2\otimes L^{-1}(P))\ge g-3$ one obtains
$h^0(L\otimes K^{-1})\ge 4$. Then, $L=K(g^3_6)$ contradicting the
assumption. If $d=2g+3$, $h^0(L\otimes K^{-1})\ge 3$ and
$L=K(g^2_5)$. In either case, this contradicts the assumptions on
$L$.

Note that if we were only assuming that $D_{g-3}$ moves in a space
of dimension $g-4$, we would obtain that $D$ has at most one fixed
point.

Assume now that $|D|$ were composed with an involution. A fixed
curve can have only a finite number of non-rational involutions.
As there is a $g-3$ dimensional family of divisors $L\otimes
K^{-1}(D_{g-3})={\mathcal O}(D)$, the involution must be rational.
So, we have
$$C\rightarrow {\bf P^1}\rightarrow {\bf P^1}$$
where the first map has degree $a$ the second degree $b$, the
composition is the map corresponding to the linear system $|D|$
and $dim |D|=b$. As $dim W^1_a\le a-2$ and $D_{g-3}$ moves in a
variety of dimension $g-3$, we would get $a-2\ge g-3$. Then, $b+1=
h^0(bg^1_a)\ge ab+1-g$. Hence, $b\le {g\over a-1}\le {g\over
g-2}<2$. Hence, $b=1$ and the map is not composed with an
involution.

This result would still hold under the assumption that $D_{g-3}$
moves in a space of dimension $g-4$.

We now prove c). We study first the case of a divisor $D_2$ of
degree two on the support of $D$. Consider the map given by the
complete linear series $|L(-D)|=|K-D_{g-3}|$. If $C$ is not
hyperelliptic, as $D_{g-3}$ is generic, this is a generic
projection of the canonical curve. Hence, it gives rise to a nodal
plane curve. Let $\bar D_2$ be the divisor corresponding to a
node, $|K-D_{g-3}-\bar D_2|\subset |K-D_{g-3}|$ the pencil cut on
$C$ by the lines through the node.

Consider the product maps
$$\begin{matrix} \mu_{\bar D_2}&:\ H^0(K-D_{g-3}-\bar D_2) \otimes H^0(K-D_2)&\rightarrow
&H^0(K^2(-D_{g-3}-D_2-\bar D_2))\\
 &\downarrow & &\downarrow\\
 \mu :\ &H^0(K(-D_{g-3}))\otimes H^0(K-D_2)&\rightarrow
&H^0(K^2(-D_{g-3}-D_2))\\
\end{matrix}$$
From the base point free pencil trick, $Ker \mu_{\bar
D_2}=H^0(D_{g-3}-D_2+\bar D_2)$.

From $degL\ge 2g+3$, one obtains $degD\ge g+2$. Hence, $|D|$ is a
linear series of dimension at least two that gives rise to a
birational map. Therefore, the monodromy of the hyperplane section
is the full symmetric group.

Hence, every pair of points of a generic divisor in $D$ imposes
the same number of conditions on $D_{g-3}+\bar D_2$. Assume
$h^0(D_{g-3}-D_2+\bar D_2)\ge 1$.  It would then follow that
$h^0(D_{g-3}+\bar D_2-D)\ge 1$ which is impossible. Therefore
$h^0(D_{g-3}+\bar D_2-D_2)=0$ and  $dim (Im \mu_{\bar
D_2})=2(g-2)=h^0(K^2(-D_{g-3}-D_2-\bar D_2))$. We prove below that
$\bar D_2$ imposes two independent conditions on $H^0(K_C-D_2)$.
Hence, it imposes at least two independent conditions on the image
of $\mu$. Then, $dim Im\mu\ge dim Im \mu_{\bar
D_2}+2=h^0(K^2(-D_{g-3}-D_2))$.

It remains to show that $\bar D_2$ imposes independent conditions
on $K-D_2$. Assume this  were not the case . Then, $\bar
D_2+D_2\in g^1_4$. If $C$ is not trigonal or bielliptic, $dim
W^1_4\le 0$. But the generic choice of $D_2$, makes this
impossible. In the bielliptic case, one of the components of
$W^1_{g-1}$ is the pull back of the set of linear series of degree
two on the elliptic curve together with sets of $g-5$ points.
Then, $\bar D_2=i(D'_2)$ with $i$ the involution on $C$ associated
to the elliptic cover. Then, $D_2+\bar D_2=g^1_4$ would imply
$D_2=D'_2$ against the condition that we checked that $D_{g-3}$
and $D$ are disjoint. Note that the condition fails for a trigonal
curve. For a generic $D_2=P+Q$, one can find a pair of points in
say $g^1_3-Q=R+S=\bar D_2$ such that $h^0(P+Q+R+S)=2$.

We now deal with the case of a divisor $D_3$ of degree three
contained in $D$. By the openess of the condition, it suffices to
prove that the map is surjective for a particular divisor
$D_{g-3}$. Choose a generic point $M\in W^2_g$. Choose a generic
effective divisor of $|M|$ and write it as $D_{g-3}+\bar D_3$.
Then, $h^0(K-D_{g-3}-\bar D_3)=2$. Consider now the commutative
diagram
$$\begin{matrix} \mu_{\bar D_3}:\ &H^0(K(-D_{g-3}-\bar D_3))\otimes H^0(K-D_3)&\rightarrow
&H^0(K^2(-D_{g-3}-D_3-\bar D_3))\\
 &\downarrow & &\downarrow\\
 \mu :\ &H^0(K(-D_{g-3}))\otimes H^0(K-D_3)&\rightarrow
&H^0(K^2(-D_{g-3}-D_3)). \\
\end{matrix}$$
 By the base point free pencil trick, the kernel of the top row is
$H^0(D_{g-3}-D_3+\bar D_3)$.

As $dim W^2_g\ge g-6$, $dim C^2_g\ge g-4$. Hence, $D_{g-3}$ moves
in a space of dimension at least $g-4$. It follows from the
remarks above that the linear series $|D|$ gives rise to a
birational map and has at most one fixed point.

 When the linear series has no fixed points,   with the same argument as in the case
of a divisor of degree two, $H^0(D_{g-3}-D_3+\bar D_3)$ has
dimension zero. Assume now that $|D|$ had a fixed point $P$. If
$D_3\subset D-P$, we deduce $dimh^0(D_{g-3}+\bar D_3-(D-P))\ge 1$
while if $D_3=D_2+P$, $dimh^0(D_{g-3}+\bar D_3-D)\ge 1$. Both are
impossible. Therefore, the image of $\mu_{\bar D_3}$ has dimension
$2(g-3)=2g-6$. We shall prove below that $\bar D_3$ imposes
independent conditions on $K-D_3$. Hence, it imposes independent
conditions on the image of the cup-product map. It follows then
that $dimIm\mu \ge 2g-6+3=h^0(K^2(-D_{g-3}-D_3-\bar D_3))$ as
required.

  As $|\bar D_3+D_{g-3}|$ is a linear series of dimension two, the monodromy
  associated to the generic hyperplane section is the whole symmetric group.
  If $\bar D_3$ does not impose independent conditions on $K-D_3$,
then any divisor of degree three  in $|\bar D_3+D_{g-3}|$ fails to
impose independent conditions on $K-D_3$. Hence, $h^0(K-D_3-\bar
D_3-D_{g-3})=g-5$ which is impossible.

\end{proof}

\end{section}

\end{document}